\documentclass[12pt]{amsart}

\usepackage{tikz}

\usepackage[latin1]{inputenc}
\usepackage{amsmath, accents} 
\usepackage{amsfonts}
\usepackage{amssymb}
\usepackage{latexsym}
\usepackage{graphicx}
\usepackage{subfigure}
\usepackage{color}
\usepackage{hyperref}
\usepackage{verbatim}
\usepackage[all]{xy}
\usepackage{graphics}
\usepackage{pdfsync}

\theoremstyle{plain}
\newtheorem{theorem}{Theorem}[section]

\newtheorem{definition}[theorem]{Definition}

\theoremstyle{remark}

\title[Generalized intersecting families]{A preliminary result for generalized intersecting families}

\author{Brian T. Chan}
\email{bchan600@gmail.com}

\date{\today}


\begin{document}

\begin{abstract} Intersecting families and blocking sets feature prominently in extremal combinatorics. We examine the following generalization of an intersecting family investigated by Hajnal, Rothschild, and others. If $s \geq 1$, $k \geq 2$, and $u \geq 1$ are integers, then say that an $s$-uniform family $\mathcal{F}$ is \emph{$(k,u)$-intersecting} if for all $A_1, A_2, \cdots, A_k \in \mathcal{F}$, $|A_i \cap A_j| \geq u$ for some $1 \leq i < j \leq k$. In this note, we investigate the following parameter. If $s$, $k$, $u$, $\ell$ are integers satisfying $s \geq 1$, $k \geq 2$, $1 \leq u \leq s$, and $2 \leq \ell < k$, then let $N^{(u)}_{k,\ell}(s)$ denote the smallest integer $r$, if it exists, such that any $(k,u)$-intersecting $s$-uniform family is the union of at most $r$ families that are $(\ell,u)$-intersecting. Using a Sunflower Lemma type argument, we prove that $N^{(u)}_{k,\ell}(s)$ always exists and that the following inequality always holds.
$$N^{(u)}_{k,\ell}(s) \; \leq \; \bigg{\lceil} \dfrac{ k - 1 }{\ell - 1} \cdot {s \choose u} \bigg{\rceil} $$
\end{abstract}
\maketitle
%
\section{Introduction} \label{sec:intro}

Intersecting families of sets and blocking sets have been of interest \cite{Jukna}. Properties of maximal intersecting families have been established by F\"{u}redi \cite{Furedi}. Moreover, Kleitman \cite{Kleitman} established an optimal upper bound for the number of sets in a union of intersecting families. Furthermore, the \emph{Erd\H{o}s-Ko-Rado Theorem} \cite{ErdosKoRado} gives an upper bound to the size of an intersecting family of $s$-element sets contained in a larger $n$-element set. Blocking sets are connected to intersecting families and have been used for Erd\H{o}s and Rado's \emph{Sunflower Lemma} and its variants and are integral to the proofs of these results \cite{ErdosRado, Furedi, HastadJuknaPudlak}. Lastly, the \emph{Kneser Conjecture} \cite{Kneser}, first proved by L\'ovasz \cite{Lovasz}, is an assertion on minimal coverings of certain families of sets by pair-wise intersecting subfamilies \cite{UnifHyp}. Subsequently, there have been extensions and generalizations of L\'ovasz's result \cite{AlonFranklLovasz, Sakaria}. \\

In this note, we consider a family of conditions for generalized intersecting families of sets previously investigated by Hajnal and Rothschild \cite{HajnalRothschild} in the context of the Erd\H{o}s-Ko-Rado theorem. Given an $s$-uniform family satisfying one of the aforementioned conditions, we are interested in the smallest number of subfamilies that exist such that the union of the subfamilies is the given family and each of the subfamilies satisfies a stronger generalized intersecting condition. We use a Sunflower Lemma type argument to establish an upper bound for this number that only depends on $s$ and parameters relating to the above intersection conditions.

\section{The upper bound}

We follow standard convention by defining, for any integer $s \geq 1$, an \emph{$s$-uniform family} to be a finite set $\mathcal{F}$ of $s$-element sets. Moreover, given an $s$-uniform family $\mathcal{F}$, we call elements of $\mathcal{F}$ a \emph{members} of $\mathcal{F}$ and subsets of $\mathcal{F}$ a \emph{subfamilies} of $\mathcal{F}$. Lastly, if $F$ is a non-empty set and if $1 \leq u \leq |F|$, then let ${F \choose u}$ denote the set of $u$-element subsets of $F$. \\

Let $s \geq 1$, $k \geq 2$, and $u \geq 1$ be integers. Then an $s$-uniform family $\mathcal{F}$ is \emph{$(k,u)$-intersecting} [cf. \cite{HajnalRothschild}] if for all $A_1, A_2, \cdots, A_k \in \mathcal{F}$, there exist integers $1 \leq i < j \leq k$ such that $|A_i \cap A_j| \geq u$. The parameter of interest in this note is as follows.

\begin{definition} For integers $s$, $k$, $u$, and $\ell$ such that $s \geq 1$, $k \geq 2$, $u \geq 1$, and $2 \leq \ell < k$. Moreover, assume that there exists a smallest positive integer $r$ such that if $\mathcal{F}$ is an $s$-uniform family that is $(k,u)$-intersecting, then there exist subfamilies $\mathcal{F}_1$, $\mathcal{F}_2$, $\cdots$, $\mathcal{F}_n$ of $\mathcal{F}$ such that $n \leq r$, $\mathcal{F} = \cup_{i=1}^n \mathcal{F}_i$, and, for all $1 \leq i \leq n$, $\mathcal{F}_i$ is $(\ell,u)$-intersecting. Then define $N^{(u)}_{k,\ell}(s) = r$. 
\end{definition}

We now prove the main result.

\begin{theorem} Let $s$, $k$, $u$, $\ell$ be integers such that $s \geq 1$, $k \geq 2$, $1 \leq u \leq s$, and $2 \leq \ell < k$. Then $N^{(u)}_{k,\ell}(s)$ exists and the following inequality holds.

$$N^{(u)}_{k,\ell}(s) \; \leq \; \bigg{\lceil} \dfrac{ k - 1 }{\ell - 1} \cdot {s \choose u} \bigg{\rceil} $$
\end{theorem}

\begin{proof} Let $s$, $k$, $u$, $\ell$ be as described in the theorem. Moreover, let $\mathcal{F}$ be an $s$-uniform family that is $(k,u)$-intersecting. By the Pigeonhole Principle, the union of $\ell - 1$ families that are $(2,u)$-intersecting is a $(\ell,u)$-intersecting family. So it is enough to prove the above theorem for the case $\ell = 2$. Because $\mathcal{F}$ is $(k,u)$-intersecting, there exist at most $k - 1$ members $F_1$, $F_2$, $\cdots$, $F_m$ of $\mathcal{F}$ such that for all $F \in \mathcal{F}$, there exists an index $1 \leq i \leq m$ such that $|F \cap F_i| \geq u$. Let $\mathfrak{X} = \cup_{i=1}^m {F_i \choose u}$. For all $X \in \mathfrak{X}$, define $\mathcal{F}_X$ to be the set of elements $F \in \mathcal{F}$ such that $X \subseteq F$. By the above, it follows that for all $F \in \mathcal{F}$, there exists an index $1 \leq i \leq m$ and a member $X \in {F_i \choose u}$ such that $X \subseteq F \cap F_i \subseteq F$. Hence, $\mathcal{F} = \bigcup_{X \in \mathfrak{X}} \mathcal{F}_X$ and we proceed as follows. For all $X \in \cup_{i=1}^m {F_i \choose u}$ and for all $F', F'' \in \mathcal{F}_X$, $X \subseteq F' \cap F''$, implying, as $|X| = u$, that $|F' \cap F''| \geq u$. So $\mathcal{F}_X$ is $(2,u)$-intersecting for all $X \in \cup_{i=1}^m {F_i \choose u}$ and the case follows since $|\cup_{i=1}^m {F_i \choose u}| \leq (k-1) \cdot {s \choose u}$.

\end{proof}

\bibliographystyle{amsplain}

\end{document}